\documentclass[submission
]{dmtcs}

\usepackage[latin1]{inputenc}
\usepackage{subfigure, amsmath, verbatim}

\usepackage[round]{natbib}

\usepackage{tikz}
\usepackage{xcolor}
\usetikzlibrary{arrows, shapes}

\newtheorem{theorem}{Theorem}
\newtheorem{proposition}[theorem]{Proposition}

\newtheorem{corollary}[theorem]{Corollary}

\newtheorem{definition}[theorem]{Definition}

\newtheorem{example}[theorem]{Example}

\newcommand{\spe}[1]{\mathtt{#1}}

\newcommand{\qbin}[2]{\begin{bmatrix}{#1}\\ {#2}\end{bmatrix}_q}

\DeclareMathOperator{\type}{type}

\DeclareMathOperator{\cone}{cone}
\DeclareMathOperator{\dcone}{dcone}
\begin{document}

\title[Triune quasisymmetric functions]{Quasisymmetric functions from combinatorial Hopf monoids and Ehrhart Theory}

\author[J.~White]{Jacob A.~White}
\address{School of Mathematical and Statistical Sciences\\
University of Texas - Rio Grande Valley\\
Edinburg, TX 78539}

\date{\today}



\maketitle

\begin{abstract}
\paragraph{Abstract.}
We investigate quasisymmetric functions coming from combinatorial Hopf monoids. We show that these invariants arise naturally in Ehrhart theory, and that some of their specializations are Hilbert functions for relative simplicial complexes. This class of complexes, called forbidden composition complexes, also forms a Hopf monoid, thus demonstrating a link between Hopf algebras, Ehrhart theory, and commutative algebra.  We also study various specializations of quasisymmetric functions.

\paragraph{Resum\'e.}
Nous \'etudions les fonctions quasisym\'etriques associ\'ees aux mono\"ides de Hopf combinatoriaux. Nous d\'emontrons que ces invariants sont des objets naturels \`a la th\'eorie de Ehrhart. De plus, certains correspondent \`a des fonctions de Hilbert associ\'ees \`a des complexes simpliciaux relatifs. Cette classe de complexes, constitue un mono\"de de Hopf, r\'ev\'elant ainsi un lien entre les alg\`ebres de Hopf, la th\'eorie de Ehrhart, et l'alg\`ebre commutative. Nous \'etudions \'egalement diverses cat\'egories de fonctions quasisym\'etriques.

\end{abstract}

\section{Introduction}

Chromatic polynomials of graphs, introduced by \cite{chromatic-polynomial} are wonderful polynomials. Their properties can be understood as coming from three different theories: 
\begin{enumerate}
\item Chromatic polynomials were shown by \cite{beck-zaslavsky} to be Ehrhart functions for inside-out polytopes. 

\item They arise as the Hilbert polynomial for the coloring ideal introduced by \cite{einar}. Moreover, this ideal is the Stanley-Reisner module for a relative simplicial complex $(\Gamma, \Delta)$.
\item The chromatic polynomial is the image of a homomorphism from the incidence Hopf algebra of graphs, first studied in \cite{Schmitt}, to the polynomial algebra.
\end{enumerate}
 Similar results have been shown by  \cite{aguiar-ardila} for the Hopf algebra of generalized Permutohedra. Thus, we can give three distinct proofs of Stanley's Reciprocity Theorem of chromatic polynomials \cite{stanley-acyclic}. Chromatic polynomials form a situation where `Ehrhart polynomial = Hilbert polynomial = polynomial coming from a Hopf algebra'. The idea of `Ehrhart = Hilbert' has been studied before by \cite{breuer-dall}. We call such polynomials \emph{triune}, because they can be studied from three different perspectives at one time.

The primary goal of this paper is to study \emph{triune} quasisymmetric functions which are Ehrhart functions, specialize to Hilbert functions, and come from combinatorial Hopf algebras. The motivation is that such invariants have three different aspects, which give them a rich structure.
Given any combinatorial Hopf monoid $\spe{H}$ with a Hopf submonoid $\spe{K}$, there is a natural quasisymmetric function $\Psi_{\spe{K}}(\spe{h})$ associated to every element $\spe{h} \in \spe{H}$. This invariant is a special case of the work of \cite{aguiar-bergeron-sottile}. In our case, the invariant can be studied from the perspective of geometric combinatorics: there is a canonical relative simplicial complex $(\Gamma_{\spe{K}, \spe{h}}, \Delta_{\spe{h}})$ associated to $\spe{h}$, with a natural geometric realization in $\mathbb{R}^I$, such that 
$\Psi_{\spe{K}}(\spe{h})$ enumerates lattice points with positive coordinates inside of the complex.
The resulting Ehrhart function is an Ehrhart quasisymmetric function as defined by \cite{breuer-klivans}. We show how principal specialization is a morphism of Hopf algebras to the ring of `Gaussian polynomial functions', and that the corresponding Ehrhart 'Gaussian polynomial' is a Hilbert function of $(\Gamma_{\spe{K}, \spe{h}}, \Delta_{\spe{h}})$ with respect to a certain bigrading. Setting $q=1$ recovers known results. 

The paper is organized as follows: we review definitions regarding the Coxeter complex of type $A$, and from Ehrhart theory.
We discuss the relationship between Ehrhart theory and Hilbert functions for relative simplicial complexes $(\Gamma, \Delta)$, where $\Delta$ is a subcomplex of the Coxeter complex, and define forbidden composition complexes. In Section 3, we review material on Hopf monoids, and define triune quasisymmetric functions, which are special cases of invariants defined by \cite{aguiar-bergeron-sottile}. In Section 4, we show that forbidden composition complexes form the terminal Hopf monoid in the category of pairs of Hopf monoids, which implies that every triune quasisymmetric function is the Ehrhart quasisymmetric function for some canonical forbidden composition complex. Thus we have a link between geometric combinatorics and combinatorial Hopf algebras that was known only in special cases. In Section 5, we discuss various specializations of quasisymmetric functions from the Hopf algebra point of view. This is motivated by the lecture notes of \cite{grinberg-reiner}, which emphasize principal specialization at $q=1$. This gives new combinatorial identities, including for Ehrhart polynomials. In the process, we discuss the notion of Gaussian polynomial function, which are linear combinations of polynomials in $q$ with $q$-binomial coefficients.

\section{Relative Composition complexes and Ehrhart Theory}

The motivation for this work comes from the study of chromatic polynomials:
\begin{enumerate}
\item In \cite{einar}, chromatic polynomials of graphs are shown to be Hilbert functions for coloring ideals, which is the Stanley-Reisner module for the relative coloring complex.
\item In \cite{beck-zaslavsky}, chromatic polynomials of graphs are shown to be Ehrhart polynomials of an inside-out polytope, which is the geometric realization of the relative coloring complex.
\end{enumerate}
Thus, the Ehrhart polynomial of the inside-out polytope of a graph is the Hilbert polynomial of coloring ideal. We give a $q$-analogue of this result for arbitrary relative composition complexes.

A \emph{set composition} is a sequence $C_1, \ldots, C_k$ of disjoint subsets of $I$ such that $\cup_{i=1}^k C_k = I$. The \emph{length} of the composition is $\ell(C) = k$. We denote set compositions with vertical bars, so $12|3$ corresponds to the set composition $\{1, 2\}, \{3\}$, and $21|3 = 12|3$. The sets $C_i$ are \emph{blocks}. Similarly, an integer composition $\alpha$ is a sequence $\alpha_1, \ldots, \alpha_k$ of positive integers whose sum is $n$. 

Given a set composition $C$, there is a natural flag of sets $F(C) := S_1 \subset S_2 \subset \cdots \subset S_{k-1} \subseteq S_k = I$, where $S_i = \cup_{j \leq i} C_i$. Similarly, given such a flag $F$, there is a set composition $C(F) := C_1, C_2, \ldots, C_k$, where $C_i = S_i \setminus S_{i-1}$. This is analogous to the classic situation for integer compositions, where there is a correspondence between integer compositions of length $\Delta$ and subsets of $[n]$ of size $k-1$. We use both notations: $S_i$ for the sets in the flag, and $C_i$ for the blocks. 
The Coxeter complex of type $A$ is the order complex on the boolean lattice $2^I \setminus I$. We let $\Sigma_I$ denote the Coxeter complex of type $A$ on the set $I$. 

\subsection{Ehrhart Quasisymmetric Function}

Let $x_1, \ldots, x_i, \ldots$ be a sequence of commuting indeterminates indexed by positive integers. A \emph{quasisymmetric function} is a power series in $x_1, \ldots$, whose terms have bounded degree, such that for any $a_1, \ldots, a_k$, and $i_1 < i_2 < \ldots < i_k$ the coefficient of $x_{i_1}^{a_1}\cdots x_{i_k}^{a_k}$ is equal to the coefficient of $x_1^{a_1}\cdots x_k^{a_k}$. A basis is given by the monomial quasisymmetric functions $M_{\alpha} = \sum_{i_1 < \cdots < i_k} \textbf{x}^{\alpha}$ where $\textbf{x}^{\alpha} = x_{i_1}^{\alpha_1}\cdots x_{i_k}^{\alpha_k}$.

Given a quasisymmetric function $Q$, and $q \in \mathbb{K} \setminus \{0\}$, and $n \in \mathbb{N}$, the \emph{principal specialization} $\mathbf{ps}(Q)(n)$ is given by $\mathbf{ps}(Q)(n) = Q(1,q,q^2,\ldots, q^{n-1}, 0, 0, 0, \ldots)$. For a fixed $q$, we view $\mathbf{ps}(Q)$ as a function from $\mathbb{N}$ to $\mathbb{K}$. When $q = 1$, we denote the specialization by $\mathbf{ps}^1(Q)(n)$. It is known that this is a polynomial function.
The \emph{stable principal specialization} is given by $\textbf{sps}(Q) = Q(1,q,q^2, \ldots)$. This gives a formal power series. However, it turns out that the coefficients $Q(n)$ of the resulting power series is a quasi-polynomial in $n$.

Given a face $F = \emptyset \subset S_1 \subset S_2 \subset \cdots \subset S_m \subset I$ of $\Sigma_I$, there is a corresponding polyhedral cone in the positive orthant $\mathbb{R}^I_{\geq 0}$. 
The cone is given by the equations $x_i < x_j$ whenever $i \in S_k, j \not\in S_k$ for some $\Delta$, and $x_i = x_j$ whenever $i \in S_k$ if and only if $j \in S_k$. For example, for the flag $\{2, 4 \} \subset \{1, 2, 4,7 \} \subset \{1, 2, 3, 4, 7, 9 \}$, we obtain the polyhedral cone given by $x_2 = x_4 < x_1 = x_7 < x_3 = x_9 < x_5 = x_8$.
Thus, for any collection $\mathcal{F}$ of faces of $\Sigma_I$, there is a collection $C(\mathcal{F})$ of open polyhedral cones in $\mathbb{R}^I_{ \geq 0}$. Given a lattice point $\textbf{a} \in \mathbb{R}^I_{>0}$, we let $\textbf{x}_{\textbf{a}} = \prod_{i \in I} x_{a_i}$ be its monomial, where the coordinates of $\textbf{a}$ are encoded in the \emph{indices}, not the exponents. 
The \emph{Ehrhart quasisymmetric function} for $C(\mathcal{F})$ is given by 
\[E_{C(\mathcal{F})} = \sum_{\mathbf{a}} \mathbf{x}_{\mathbf{a}} \]
where the sum is over all lattice points which lie in some cone of $C(\mathcal{F})$.
Since $E_{C(\mathcal{F})} = \sum_{F \in \mathcal{F}} M_{\type(C(F))}$, this is a quasisymmetric function, first appearing in the work of  \cite{breuer}.

We mention specializations of $E_{C(\mathcal{F})}$, and their combinatorial interpretations. First, $\textbf{ps}^1(E_{C(\mathcal{F})})(n+1)$ is the number of lattice points in $C(\mathcal{F}) \cap [0,n]^{|I|}$. Also, $[q^m]\textbf{sps}(E_{C(\mathcal{F})})$ is the number of lattices points in $C(\mathcal{F}) \cap \Gamma_m$, where $\Gamma_m$ is the simplex given by the equation $\sum_{i \in I} a_i = m$. Finally, $[q^m] \textbf{ps}(E_{C(\mathcal{F})})(n+1)$ is the number of lattice points in $C(\mathcal{F}) \cap \Gamma_m \cap [0,n]^{|I|}$.

\subsection{Relative Composition complexes}

We define Stanley-Reisner modules for relative simplicial complexes, and introduce relative composition complexes, which have a natural geometric realization as open polyhedral cones. We show that specializations of the Ehrhart quasisymmetric function for relative composition complexes give the Hilbert function of the Stanley-Reisner module.

A \emph{relative simplicial complex} is a pair $(\Gamma, \Delta)$ where $\Gamma \subseteq \Delta$, and $\Delta$ is a simplicial complex. Given $\Delta$ with vertices $S$, we let $\mathbb{C}[S]$ be the polynomial ring with indeterminates $s_1, \ldots, s_k$, the vertices of $S$. The \emph{Stanley-Reisner ideal} for $\Delta$ is generated by $\langle \sigma \subseteq S : \sigma \not\in \Delta \rangle$, and the \emph{Stanley-Reisner module} for $(\Gamma,\Delta)$ is $I_{\Gamma} / I_{\Delta}$. The module is graded by total degree, and its Hilbert function $H(\Gamma, \Delta)(n)$ is the number of monomials of degree $n$ in the module. It is known that the Hilbert function is in fact a polynomial: details can be found in \cite{stanley-algebra}. 

We also extend the definition of double coning over a simplicial complex. Coning over a non-void complex $\Delta$ consists of adding new vertex $x$, and adding new faces $\sigma \cup \{x\}$ for all $\sigma \in \Delta$. If $\Delta = \emptyset$, we let $\cone(\Delta) = \emptyset $.
For $(\Gamma, \Delta)$, $\cone(\Gamma, \Delta) = (\cone(\Gamma), \cone(\Delta))$. Finally, the double cone is defined by $\dcone(\Gamma, \Delta) = \cone(\cone(\Gamma, \Delta))$.

Now we define the simplicial complexes that are of interest to us. A \emph{relative composition complex} is a relative complex $(\Gamma, \Delta)$ where $\Delta \subseteq \Sigma_I$. 
Given a relative composition complex $(\Gamma, \Delta)$, a composition $C$ of $\Gamma$, and a block $B$ of $C$, $B$ is \emph{forbidden} if every composition of $\Delta$ that refines $C$ does not contain $B$ as block. $(\Gamma, \Delta)$ is a \emph{forbidden complex} if every composition of $\Gamma$ either has a forbidden block, or is a facet of $\Delta$. While the definition seems unusual, we will see that forbidden composition complexes arise naturally in the study of Hopf monoids in species.

\begin{example}
\label{ex:comp}
Here are some examples.
\begin{enumerate}
\item  Let $\Gamma_{k,I}$ consist of all set compositions that have at least one block of size $\geq k$. Then points in $\Gamma_{k,I}$ consist of points in $\mathbb{R}^I$ that have at least $k$ equal coordinates. This arises in the study of the $k$-equal problem. We see that $(\Gamma_{k,I}, \Sigma_I)$ is a forbidden composition complex.

\item Let $I = \{a,b,c,d \}$, and let $\Delta$ be the complex with facets which correspond to the permutations $abcd, abdc, adbc, adcb, dabc, dacb, dcab, dcba$. Let $\Gamma$ be the subcomplex with facets $cd|b|a, cd|a|b, ab|c|d$, $ab|d|c, a|bc|d, a|cd|b, d|ab|c$, and $d|bc|a$. Then $(\Gamma, \Delta)$ is a forbidden composition complex. The Ehrhart quasisymmetric function is $8M_{1111} + 4M_{112} + 2M_{121} + 2M_{211} + M_{22}$.

\item Let $\Gamma$ be the simplex corresponding to $12|34|56$. Then $(\Gamma, \Sigma_{[6]})$ is a relative composition complex, but it is not a forbidden complex, because $1|2|34|56$, $12|3|4|56$ and $|12|34|5|6$ are all faces of $\Sigma_{[6]}$. Thus, $12|34|56$ has no forbidden blocks.

\end{enumerate}

\end{example}

Forbidden composition complexes generalize coloring complexes.
Given a graph $\spe{g}$, let $\Gamma_{\spe{g}}$ denote the collection of set compositions $C$ for which some block contains an edge of $\spe{g}$. This is the coloring complex introduced by \cite{einar}. We let $(\Gamma_{\spe{g}}, \Sigma_I)$ be the relative coloring complex. The Stanley-Reisner module for the double cone over $(\Gamma_{\spe{g}}, \Sigma_I)$ is the coloring ideal.
Our relative coloring complex is thus an example of a forbidden composition complex.

In the case of relative composition complexes, the polynomial ring associated to $\Delta$ has indeterminates given by all subsets $S \subseteq I$. We define the bidegree of $S$ to be $(|S|, 1)$. In this case, $H(\Gamma, \Delta)(m,n)$ is the function which counts the number of monomials of degree $(m,n)$. 
\begin{theorem}
Let $(\Gamma, \Delta)$ be a relative composition complex. Then $H(\dcone(\Gamma, \Delta))(n) = \textbf{ps}^1(E_{C_{\Delta \setminus \Gamma}})(n)$. Similarly, letting $H(\dcone(\Gamma, \Delta))(q,n) = \sum_{m \geq 0} H(\dcone(\Gamma, \Delta))(m,n) q^m$, we have \\
 $q^{n|I|}H(\dcone(\Gamma, \Delta))(q^{-1},n) = \textbf{ps}(E_{C(\Delta \setminus \Gamma)})(n)$.
\end{theorem}
Our first result follows from work of \cite{breuer-klivans}. However, in their setting there is no natural Stanley-Reisner module. The second result is similar to work of \cite{breuer-dall}.

\section{Hopf monoids and Characters}
In this section, we dicuss combinatorial Hopf monoids, their characters, and their quasisymmetric functions. Hopf monoids are a generalization of graphs, posets and matroids. The idea is that we have some notion of combinatorial structure, called a species, as introduced by \cite{joyal}. Moreover, we have rules for combining and decomposing these structures in a coherent way. Hopf monoids in species were originally introduced in \cite{aguiar-mahajan-1}, although the variation we discuss here can be found in \cite{aguiar-mahajan-2}. Hopf monoids allow us to define a whole class of quasisymmetric functions, and prove identities relating quasisymmetric functions in the same class, such as the class of chromatic symmetric functions of graphs.

\subsection{Hopf monoids in species}

\begin{definition}
A \emph{species} is an endofunctor $\spe{F}: \mbox{Set} \to \mbox{Set}$ on the category of finite sets with bijections.
For each finite set $I$, $\spe{F}_I$ is a finite set, and for every bijection $\sigma:I \to J$ between finite sets, there is a bijection $\spe{F}_{\sigma}: \spe{F}_I \to \spe{F}_J$, such that $\spe{F}_{\sigma \circ \tau} = \spe{F}_{\sigma} \circ \spe{F}_{\tau}$ for every pair $\sigma: I \to J$, $\tau: K \to I$.
 It is \emph{connected} if $|\spe{F}_{\emptyset}| = 1$. All species in this paper are connected, and $1_{\spe{F}}$ denotes the only element of $\spe{F}_{\emptyset}$.
\end{definition}
The exponential generating function for $\spe{F}$ is $\spe{F}(x) = \sum_{n \in \mathbb{N}} |\spe{F}_{[n]}| \frac{x^n}{n!}$.

\begin{example}
We list various examples of species.
\begin{enumerate}
\item The graph species $\spe{G}$: the set $\spe{G}_I$ consists of all graphs with vertex set $I$. Given $\sigma:I \to J$, and $\spe{g} \in \spe{G}_I$, $\spe{G}_{\sigma}(\spe{g}) = \spe{h}$ is the graph on vertex set $J$ where $i \sim j$ in $\spe{h}$ if and only if $\sigma^{-1}(i) \sim \sigma^{-1}(j)$ in $\spe{g}$. Then $\spe{G}(x) = \sum_{n \geq 0} 2^{\binom{n}{2}} \frac{x^n}{n!} = 1 + x + 2 \frac{x^2}{2} + 8 \frac{x^2}{6} + \cdots$
\item The poset species $\spe{P}$: the set $\spe{P}_I$ consists of all partial orders on $I$.  Given $\sigma:I \to J$, and $\spe{p} \in \spe{P}_I$, $\spe{G}_{\sigma}(\spe{p}) = \spe{q}$ is the partial order on $J$ where $i \leq_{\spe{q}} j$ if and only if $\sigma^{-1}(i) \leq_{\spe{p}} \sigma^{-1}(j)$. $\spe{P}(x) = 1 + x + 3\frac{x^2}{2} + 19\frac{x^3}{6} + \cdots$. 

\item The matroid species, whose structures $\spe{M}_I$ consist of all matroids on $I$. Then $\spe{M}(x) = 1 + 2x + 5 \frac{x^2}{2} + 16 \frac{x^3}{6} + \cdots$
\item The species $\spe{R}$: the set $\spe{R}_I$ consists of all relative composition complexes $(\Gamma, \Delta)$ where $\Delta \subseteq \Sigma_I$. 
\end{enumerate}

\end{example}

\begin{definition}
A \emph{monoid} is a species $\spe{F}$, equipped with associative multiplication maps $\mu_{S,T}: \spe{F}_S \times \spe{F}_T \to \spe{F}_{S \sqcup T}$ for every pair $S,T$ of finite sets, where $S \sqcup T$ denotes disjoint union. We denote the product of $\spe{f} \in \spe{F}_S$, $\spe{g} \in \spe{F}_T$ by $\spe{f} \cdot \spe{g}$. Associativity means that $(\spe{f} \cdot \spe{g}) \cdot \spe{h} = \spe{f} \cdot (\spe{g} \cdot \spe{h})$ whenever the multiplication is defined. Moreover, $1_{\spe{F}} \cdot \spe{f} = \spe{f} = \spe{f} \cdot 1_{\spe{F}}$.
\end{definition}

\begin{example}
We list various monoid operations.
\begin{enumerate}
\item The graph species $\spe{G}$ is a monoid. Given two graphs $\spe{g}$ and $\spe{h}$ with disjoint vertex sets, $\spe{g} \cdot \spe{h}$ is their disjoint union: the graph with edges $i \sim j$ if and only if $i,j \in V(\spe{g})$ and $i \sim j$ in $\spe{g}$, or $i,j \in V(\spe{h})$, and $i \sim j$ in $\spe{h}$. 
\item The poset species $\spe{P}$ is a monoid. The product is also given by disjoint union of partial orders.
\item The matroid species $\spe{M}$ is a monoid. The product is the direct sum.
\item The relative composition complex species $\spe{R}$ is a monoid. Given $(\Gamma, \Delta) \in \spe{R}_I$ and $(\Gamma', \Delta') \in \spe{R}_J$, we let $\Delta \cdot \Delta'$ be the set of all quasi-shuffles of set compositions $C, C'$, where $C \in \Delta$ and $C' \in \Delta''$. Here is an example: the quasi-shuffles of $1| 2$ and $a|b$ are: $1|2|a|b, 1|2a|b,  1|a|2|b, 1|a|2b$, $1a|2|b, 1a|2b, a|1|2|b, a|1|2b$, $a|1|b|2,  a|1b|2, a|b|1|2, 1|a|b|2$, and $1a|b|2$. The subcomplex $\Gamma \cdot \Gamma'$ consists of all quasi-shuffles of all compositions $C, C'$ where $C \in \Gamma, C' \in \Delta'$ or $C \in \Delta$ and $C' \in \Gamma'$. The product is then $(\Gamma, \Delta) \cdot (\Gamma', \Delta') = (\Gamma \cdot \Gamma', \Delta \cdot \Delta')$.
\item The species $\spe{\Phi}$ of forbidden composition complexes is a submonoid of $\spe{R}$.
\end{enumerate}
\end{example}

\begin{definition}
A \emph{combinatorial Hopf monoid} (in species) is a monoid $\spe{F}$ such that, for every $S \subseteq I$, there are \emph{partial} functions restriction $\spe{F}|_S: \spe{F}_I \to \spe{F}_S$ and contraction $\spe{F}/S: \spe{F}_I \to \spe{F}_{I \setminus S}$, subject to:
\begin{enumerate}
\item For any $T \subseteq S \subseteq I$, $(\spe{f}|_S)|_T = \spe{f}|_T$.
\item For any $T \subseteq S \subseteq I$, $(\spe{f}/T)/S = \spe{f}/S$.
\item For any $T \subseteq S \subseteq I$, $(\spe{f}|_S)/T = (\spe{f}/T)|_S$.
\item $(\spe{f} \cdot \spe{g})|_S = \spe{f}|_{S \cap A} \cdot \spe{g}|_{S \cap B}$.
\item $(\spe{f} \cdot \spe{g})/S = \spe{f}/(S \cap A) \cdot \spe{g}/(S \cap B)$.
\end{enumerate}
 
\end{definition}
We are working with partial functions, so if one side of the equation is undefined, then so is the other side.

\begin{example}
We list various examples of Hopf monoids.
\begin{enumerate}
\item The graph species $\spe{G}$ is a combinatorial Hopf monoid. The restriction map $\spe{g}|_S$ consists of the graph on $S$ with edges $i \sim j$ if and only if $i \sim j$ in $\spe{g}$. In this case, we define $\spe{g} / S = \spe{g}|_{I \setminus S}$. 
\item The poset species $\spe{P}$ is a combinatorial Hopf monoid. Given $\spe{p}$, and $S \subseteq I$, we let $\spe{p}[I]$ denote the induced subposet. Then $\spe{p}|_S = \spe{p}[S]$ provided $\spe{p}[S]$ is an order ideal of $\spe{p}$. Similarly, if $\spe{p}[S]$ is a lower order ideal, we let $\spe{p}/S = \spe{p}[I \setminus S]$.
\item The matroid species $\spe{M}$ is a combinatorial Hopf monoid, with the restriction and contraction.
\item The species $\spe{\Phi}$ of forbidden composition complexes is a Hopf monoid. We do not give the definition of the coproduct, as it is technical. It involves the notion of deconcatenation.
\end{enumerate}
\end{example}

\begin{figure}[htbp]
\begin{center}
\begin{tabular}{l|c|c|c}
Species & Example $x$ & $x|_S$ & $x/S$ \\

$\spe{G}_I$ &

\begin{tikzpicture}

\node[draw=black, fill=magenta, circle] (x) at (0,0) {x};
\node[draw=black, fill=magenta, circle] (y)  at (2,0) {y};
\node[draw=black, circle] (z)  at (4,0) {z};

\draw[black, line width = 3pt]  (x) -- (y);
\draw[black, line width = 3pt]  (y) -- (z);
\end{tikzpicture}

&

\begin{tikzpicture}

\node[draw=black, fill=magenta, circle] (x) at (0,0) {x};
\node[draw=black, fill=magenta, circle] (y)  at (2,0) {y};

\draw[black, line width = 3pt]  (x) -- (y);

\end{tikzpicture}

&

\begin{tikzpicture}

\node[draw=black, circle] (z)  at (4,0) {z};

\end{tikzpicture}

\\

$\spe{P}_I$ &

\begin{tikzpicture}[>=latex]

\node[draw=black, fill=magenta, circle] (x) at (0,0) {x};
\node[draw=black, circle] (y)  at (2,2) {y};
\node[draw=black, circle] (z)  at (4,0) {z};

\draw[black, line width = 3pt, ->] (y) -- (x);
\draw[black, line width = 3pt, ->] (y) -- (z);
\end{tikzpicture}

&

\begin{tikzpicture}[>=latex]

\node[draw=black, fill=magenta, circle] (x) at (0,0) {x};

\end{tikzpicture}

&

\begin{tikzpicture}[>=latex]

\node[draw=black, circle] (y)  at (0,2) {y};
\node[draw=black, circle] (z)  at (0,0) {z};

\draw[black, line width = 3pt, ->] (y) -- (z);
\end{tikzpicture}

\\

$\spe{P}_I$ & 

\begin{tikzpicture}[>=latex]

\node[draw=black, circle] (x) at (0,0) {x};
\node[draw=black, fill=magenta, circle] (y)  at (2,2) {y};
\node[draw=black, circle] (z)  at (4,0) {z};

\draw[black, line width = 3pt, ->] (y) -- (x);
\draw[black, line width = 3pt, ->] (y) -- (z);
\end{tikzpicture}

& 
undefined: $x|_S$ is not an order ideal.
& undefined

\end{tabular}
\end{center}
\caption{Examples of restriction and quotient. Shaded vertices are elements of $S$.}

\end{figure}
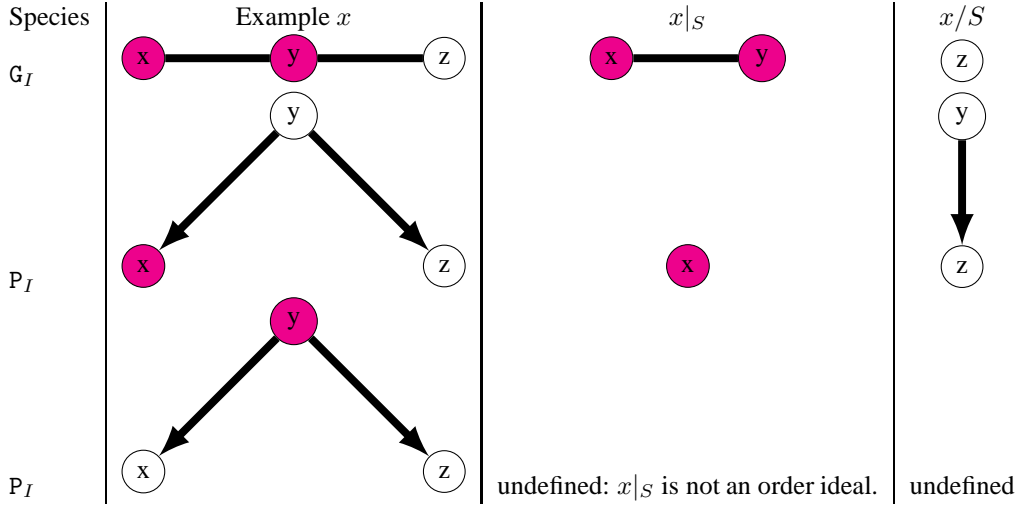

\subsection{Characters and Inversion}
Now we discuss characters of Hopf monoids.
\begin{definition}
Given a Hopf monoid $\spe{H}$, and a field $\mathbb{K}$, a character is a multiplicative function $\varphi: \spe{H} \to \mathbb{K}$. For every finite set $I$, there is a map $\varphi_I: \spe{H}_I \to \mathbb{K}$, natural in $I$, such that, for all $I = S \sqcup T$, $\spe{h}_S \in \spe{H}_S, \spe{h}_T \in \spe{H}_T$, we have $\varphi_S(\spe{h}_S)\varphi_T(\spe{h}_T) = \varphi_I(\spe{h}_S \cdot \spe{h}_T)$. The character is \emph{connected} if $\varphi_{\emptyset}(\spe{h}_{\emptyset}) = 1$.
\end{definition}

\begin{example}
One example is the character given by $\varphi_I(\spe{h}) = 1$ for all $I$, $\spe{h} \in \spe{H}_I$. This is the \emph{zeta} character.

 Let $\spe{G}$ be the Hopf monoid of graphs. Given a graph $\spe{g}$, let \[\varphi(\spe{g}) =  \left\{\begin{array}{cc} 1 & \spe{g} \mbox{ has no edges} \\ 0 & \mbox{ otherwise } \end{array}\right. \]
\end{example}

The set  $\chi(\spe{H})$ of connected characters on $\spe{H}$ is a group, with multiplication given by: 
 \[(\varphi \ast \psi)_I(\spe{h}) = \sum_{S \subseteq I} \varphi_S(\spe{h}|_S) \psi_{I \setminus S}(\spe{h}/S) \]
for $\varphi, \psi \in \chi(\spe{H})$,
where the right hand side is $0$ for any $S$ where $h|_S$ or $h/S$ is undefined.

The inverse of a character $\varphi$ is defined recursively:
\begin{enumerate}
\item $\varphi^{-1}_{\emptyset} = \varphi_{\emptyset}$
\item For $\spe{h} \in \spe{H}_I$, $\varphi^{-1}(\spe{h}) = - \sum\limits_{S \subset I} \varphi^{-1}_S(\spe{h}|_S) \varphi_{I \setminus S} (\spe{h}/S)$
\end{enumerate}

We discuss characters coming from Hopf submonoids $\spe{K} \subseteq \spe{H}$. 
A \emph{Hopf submonoid} $\spe{K}$  is a subspecies, meaning that $\spe{K}_I \subseteq \spe{H}_I$ for all $I$. Moreover, the product, restriction, and contraction of elements of $\spe{K}$ remain in $\spe{K}$.
Given a submonoid $\spe{K} \subset \spe{H}$, there is a character $\varphi_{\spe{K}}: \spe{H} \to \mathbb{K}$ given by: 
\[\varphi_{\spe{K}}(\spe{h}) = \left\{\begin{array}{cc} 1 & \spe{h} \in \spe{K}_I \\ 0 & \mbox{ otherwise } \end{array}\right. \]

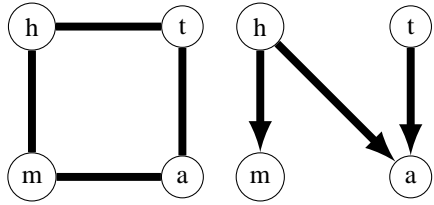
\begin{figure}[htpb]
\label{fig:graphposet}
\begin{center}
\begin{tabular}{cc}
\begin{tikzpicture}

\node[draw=black, circle] (x) at (0,0) {m};
\node[draw=black, circle] (y)  at (2,0) {a};
\node[draw=black, circle] (z)  at (2,2) {t};
\node[draw=black, circle] (w) at (0,2) {h};
\draw[black, line width = 3pt]  (x) -- (y);
\draw[black, line width = 3pt]  (y) -- (z);
\draw[black, line width = 3pt]  (z) -- (w);
\draw[black, line width = 3pt]  (w) -- (x);
\end{tikzpicture}

&
\begin{tikzpicture}[>=latex]

\node[draw=black, circle] (x) at (0,0) {m};
\node[draw=black, circle] (y)  at (2,0) {a};
\node[draw=black, circle] (z)  at (2,2) {t};
\node[draw=black, circle] (w) at (0,2) {h};
\draw[black, line width = 3pt, ->]  (w) -- (y);
\draw[black, line width = 3pt, ->]  (z) -- (y);

\draw[black, line width = 3pt, ->]  (w) -- (x);
\end{tikzpicture}
\end{tabular}
\end{center}
\caption{Example graph and poset}
\end{figure}

In the case of graphs, one Hopf submonoid is the species of edgeless graphs. In the case of posets, there is the Hopf submonoid of antichains. There is a Hopf monoid of generalized Permutohedra, and the character \cite{aguiar-ardila} study also comes from a Hopf submonoid. Finally, the Hopf monoid of composition complexes is a Hopf submonoid of $\spe{\Phi}$. In each of these cases, we obtain a character.

\subsection{The quasisymmetric function associated to a character}

We recall the quasisymmetric function associated to a character on a Hopf monoid $\spe{H}$. Given a set composition $C$ of $I$, and $\spe{h} \in \spe{H}_I$, define $\varphi_C(h_1, \ldots, h_k) = \prod_{i=1}^k \varphi(h_i)$,  $h_i = (h|_{S_i}) / S_{i-1}$, where $S_i \in F(C)$.

\begin{definition}
Given a combinatorial Hopf monoid $\spe{H}$, a character $\varphi$, a finite set $I$, and $\spe{h} \in \spe{H}_I$, define \[\Psi_{\varphi}(\spe{h}) = \sum_{C \models I} \left(\prod_{i=1}^{\ell(C)} \varphi(\spe{h}|_i)\right) M_{\type(C)}.\]
\end{definition}
 Given a combinatorial Hopf monoid, the vector space generated by the equivalence classes of $\spe{H}$-structures forms a combinatorial Hopf algebra, which appears in \cite{aguiar-mahajan-1}. Moreover, by work of \cite{aguiar-bergeron-sottile}, there is a unique morphism from this Hopf algebra to $QSym$. Our definition $\Psi_{\varphi}(\spe{H})$ is the resulting map.

There is a description for $\Psi$ in terms of colorings. Given a coloring $f: I \to \mathbb{N}$, and $i \in \mathbb{N}$, we let $\spe{h}|_i = \spe{h}|_{f^{-1}([i])} / f^{-1}([i-1])$ be the $i$th minor of $\spe{h}$ under $f$.
\begin{theorem}
Let $\spe{H}$ be a combinatorial Hopf monoid, with a character $\varphi: \spe{H} \to \spe{E}$. Fix a finite set $I$, and $\spe{h} \in \spe{H}_I$. 
Then  \[\Psi_{\varphi}(\spe{h}) = \sum_{f: I \to \mathbb{N}} \varphi_f(\spe{h}) \mathbf{x}_f \]
where $\varphi_f(\spe{h}) = \prod_{i \in \mathbb{N}} \varphi(\spe{h}|_i)$, which is well-defined.

\end{theorem}

For a coloring of a graph $\spe{g}$, , the $i$th minor is the induced subgraph on the $i$th color class, so $\varphi(\spe{g}|_i) = 1$ if and only if the $i$th color class is an independent set. Thus our quasisymmetric function enumerates proper colorings, giving the chromatic symmetric function introduced by \cite{stanley-coloring-1}. For example, for the graph in figure \ref{fig:graphposet} the resulting chromatic symmetric function is $24M_{1111} + 4M_{211} + 4M_{121} + 4M_{112} + 2M_{22}$.
For posets, $\varphi_f(\spe{p}) = 1$ if and only if $f: I \to \mathbb{N}$ is a strictly order preserving map, which is the quasisymmetric function for strict $P$-partitions considered by \cite{stanley-poset}. For example, for the poset $p$ in figure \ref{fig:graphposet}, the quasisymmetric function is given by $5M_{1111} + 2M_{211} + M_{121} + 2M_{112} + M_{22}$.

\begin{theorem}
Let $\spe{\Phi}$ be the Hopf monoid of forbidden composition complexes. Then for all $(\Gamma, \Delta) \in \spe{\Phi}_I$, we have $\Psi_{\varphi}(\Gamma, \Delta) = E_{C(\Delta \setminus \Gamma)}$.
\end{theorem}

\section{Forbidden composition complexes}
We show that, for any pair $\spe{K} \subseteq \spe{H}$ of Hopf monoids, and any element $\spe{h} \in \spe{H}_I$, there is a forbidden composition complex $(\Gamma_{\spe{K}, \spe{h}}, \Delta_{\spe{h}})$ whose Ehrhart quasisymmetric function is $\Psi_{\varphi_{\spe{K}}}(\spe{h})$. 
Let $\spe{H}$ be a combinatorial Hopf monoid, $\spe{h} \in \spe{H}_I$, and define $\Delta_{\spe{h}} \subseteq \Sigma$ to be those faces $F$ such that $h_i$ is defined for all $S_i \in F$. This defines a morphism of Hopf monoids $\Delta: \spe{H} \to \spe{C}$, the species of composition complexes.

Let $\spe{K} \subseteq \spe{H}$ be a Hopf submonoid,
and define $\Delta_{\spe{K}, \spe{h}}$ to consist faces $F \in \Delta_{\spe{h}}$ such that some minor $h|_i \not\in \spe{K}_{S_i - S_{i-1}}$. Since $\spe{K}$ is a Hopf submonoid, $\Gamma_{\spe{K}, \spe{h}}$ is a simplicial complex. Moreover, $(\Gamma_{\spe{K}, \spe{h}}, \Delta_{\spe{h}})$ is a forbidden composition complex, and the map $\Gamma_{\spe{K}}: \spe{H} \to \spe{\Phi}$ defined by $\Gamma(\spe{h}) = (\Gamma_{\spe{K}, \spe{h}}, \Delta_{\spe{H}})$ is a morphism of Hopf monoids.
However, even more is true: the quasisymmetric function $\Psi_{\varphi_{\spe{K}}}(\spe{h})$ is the Ehrhart quasisymmetric function for $\Gamma(\spe{h})$:
\begin{theorem}
Given a set $I$, let $\spe{\Phi}_I$ denote the set of all forbidden composition complexes on $I$, and $\spe{C}_I$ denote the set of all composition complexes on $I$.
\begin{enumerate}
\item Given any combinatorial Hopf monoid $\spe{H}$, with $\spe{K} \subseteq \spe{H}$, there exists unique morphisms of combinatorial Hopf monoids $\Delta: \spe{K} \to \spe{C}$, $\Gamma_{\spe{K}}: \spe{H} \to \spe{\Phi}$ such that $\iota \Delta = \Gamma_{\spe{K}} \iota$, where $\iota$ denotes inclusion maps $\spe{K} \subseteq \spe{H}$, and $\spe{C} \subseteq \spe{\Phi}$.
\item Under this map, $ \Psi_{\varphi_{\spe{K}}}(\spe{h}) = E(C(\Gamma(\spe{h})))$.
\end{enumerate}
\end{theorem}
For a graph $\spe{g}$, $\Delta_{\spe{g}}$ is the relative coloring complex. Given a poset $\spe{p}$, let $C(\spe{p})$ be the polyhedral cone in $\mathbb{R}^I$ bounded by equations $x_i \leq x_j$, for all $i \leq j$ in $\spe{p}$. Then $\Delta_{\spe{p}}$ consists of all cones in the Coxeter arrangement which lie in $C(\spe{p})$. Similarly, $\Gamma_{\spe{K}, \spe{p}}$ consists of the cones which lie on the boundary of $C(\spe{p})$.

\section{Specializations}

We discuss specializations of quasisymmetric functions, and interpretations of $\Psi$ under specialization. We show combinatorial identities relating quasisymmetric functions for various elements of the same combinatorial Hopf monoid.
It is known that $\textbf{ps}^1$ is a Hopf algebra homomorphism from $QSym$ to $\mathbb{K}[x]$. We show that $\mathbf{ps}$ is a morphism of Hopf algebras in general. The image of $\mathbf{ps}$ is the ring of Gaussian polynomial functions, which are $q$-analogues of polynomials. We also study the stable principal specialization $\textbf{psp}$. While this section primarily emphasizes the Hopf algebra perspective, many of the results are of combinatorial interest.

\subsection{Gaussian polynomials and principal specialization}

Clearly, $\mathbf{ps}(Q): \mathbb{N} \to \mathbb{C}$ is a polynomial function when $q=1$. This leads to the question of what type of function we get for general $q$. For now, assume that we are working over $\mathbb{C}(q)$.

 For any integer $m$, define $D_m(f):\mathbb{N} \to \mathbb{C}(q)$ by $D_m(f)(n) = f(n+1) - q^mf(n)$, and $D^m(f) =D_m \circ D^{m-1}(f)$. A function $f$ is a \emph{Gaussian polynomial function} of degree at most $d$ if $D^{d+1}(f) = 0$. We recovering the classical definitions when $q = 1$. The terminology comes from the fact that $q$-binomial coefficients are sometimes called Gaussian polynomials, and all Gaussian polynomial functions can be expressed as linear combinations of $q$-binomial coefficients.
Consider a Gaussian polynomial function of degree $m$. Then we can define $f(-n) = q^{-m}(f(-n+1) - D_m(f)(-n))$, for $n > 0$. Thus Gaussian polynomials are functions from $\mathbb{Z} \to \mathbb{C}$. 

\begin{theorem}
The algebra of Gaussian polynomials, $G$, is a Hopf algebra, with basis given by $[x]^n$, $n \in \mathbb{N}$. The unit is $1$, and multiplication is given by $[x]^k \cdot [x]^m = [x]^{k+m}$.  The comultiplication sends $[x]$ to  $[x] \otimes 1 + q^{x} \otimes [x]$, and the antipode is generated by $S([x]) = [-x]$.

Moreover, $\textbf{ps}: QSym \to G$ is a morphism of Hopf algebras, and $G$ is graded as an algebra, but not as a coalgebra.
\end{theorem}

Let us consider some examples. For the graph $\spe{g}$ in figure \ref{fig:graphposet}, the resulting Gaussian chromatic polynomial is $14q^6 \qbin{n}{4} + (2q^5+4q^4+2q^3)\qbin{n+1}{4} + 2q^2\qbin{n+2}{4}$.
For the poset in figure \ref{fig:graphposet}, the resulting Gaussian polynomial is $q^6 \qbin{n}{4} + (q^5+q^4+q^3)\qbin{n+1}{4} + q^2\qbin{n+2}{4}$. Finally, for the forbidden composition complex in Example \ref{ex:comp}, part 2, the resulting Gaussian polynomial is
$2q^6 \qbin{n}{4} + (3q^5+q^4+q^3)\qbin{n+1}{4} + q^2\qbin{n+2}{4}$.

\begin{definition}
Given a combinatorial Hopf monoid $\spe{H}$ with a character $\varphi$, and $\spe{h} \in \spe{H}_I$, define $P_{\varphi}(\spe{h},q,n) = \textbf{ps}(\Psi_{\varphi}(\spe{h}))$. 
This is the polynomial of $\spe{h}$ associated to $\varphi$.
Alternatively, $P_{\varphi}(\spe{h},q,n) = \sum_{f:I \to [n]} \varphi_f(\spe{h}) q^{w(f)}$, where $w(f) = \sum_{i \in I} (f(i) - 1)$.
\end{definition}
For a poset $\spe{p}$, $[q^n]P_{\varphi}(\spe{p},q,m)$ is the number of $\spe{p}^{\ast}$-partitions of $n$ with part size at most $m$.
\begin{proposition}
\label{prop:additive}
Let $\spe{H}$ be a combinatorial Hopf monoid with character $\varphi$, and let $\spe{h} \in \spe{H}_I$, $\spe{k} \in \spe{H}_J$, where $I$ and $J$ are disjoint sets.
Then the following identities hold:
\begin{enumerate}
\item $P_{\varphi}(\spe{h} \cdot \spe{k}, q, n) = P_{\varphi}(\spe{h}, q, n) \cdot P_{\varphi}(\spe{k},q,n)$.
\item for any $n,m \in \mathbb{N}$, $P_{\varphi}(\spe{h}, q,m+n) = q^{m|I|}\sum_{S \subseteq I}q^{-m|S|} P_{\varphi}(\spe{h}|_S, q, m) \cdot P_{\varphi}(\spe{h}/S, q, n)$.
\item  $P_{\varphi}(\spe{h}, q,-n) = q^{|I|}P_{\varphi^{-1}}(\spe{h}, q^{-1}, n)$.
\end{enumerate}
\end{proposition}
The last identity is a reciprocity result. On the left-hand side, we are counting negative colors, so we expect to have negative powers of $q$. When $q=1$, these identities are already known for graphs and posets. Also, combinatorial reciprocity for $P_{\varphi}(\spe{p}, q, n)$ is also due to \cite{stanley-poset}.

When $\varphi = \varphi_{\spe{K}}$ for some Hopf submonoid $\spe{K} \subseteq \spe{H}$, then the fact that $P_{\varphi_{\spe{K}}}(\spe{h}, 1, n)$ is an Ehrhart function, and a Hilbert function, allows us to conclude new results regarding $\varphi^{-1}$.
\begin{corollary}
Let $\spe{K} \subseteq \spe{H}$ be an inclusion of Hopf monoids, and let $\spe{h} \in \spe{H}_I$. Then $(-1)^{|I|}P_{\varphi^{-1}_{\spe{K}}}(\spe{h}, 1, n) = \sum_{\mathbf{a} \in \mathbb{Z}^I \cap (0,n]^I} w(\mathbf{a})$, where $w(\mathbf{a})$ is the number of cones $C$ of $(\Gamma_{\spe{K}, \spe{h}}, \Delta_{\spe{h}})$ such that $\mathbf{a} \in \bar{C}$.

Moreover, $\varphi^{-1}(\spe{h}) = \chi(\Gamma_{\spe{K}, \spe{h}}, \Delta_{\spe{h}}) = \sum_{\sigma \in \Delta_{\spe{h}} \setminus \Gamma_{\spe{K}, \spe{h}}} (-1)^{|\sigma|}$, the Euler characteristic.
\end{corollary}

\subsection{The stable principal specialization}

We define $Q_{\varphi}(\spe{h}, q) = \textbf{sps}(\Psi_{\varphi}(\spe{h}))$. In particular, $Q_{\varphi}(\spe{h}, q) = \sum_{f:I \to \mathbb{N}} \varphi_f(\spe{h}) q^{w(f)}$.  For posets, $Q_{\varphi}(\spe{p},q)$ is the generating function for strict $\spe{p}^{\ast}$-partitions.

\begin{theorem}

Let $\spe{H}$ be a combinatorial Hopf monoid with character $\varphi$, and let $\spe{h} \in \spe{H}_I$, $\spe{k} \in \spe{H}_J$, where $I$ and $J$ are disjoint sets.
Then the following identities hold:
\begin{enumerate}
\item $Q_{\varphi}(\spe{h} \cdot \spe{k},q) = Q_{\varphi}(\spe{h},  q) \cdot Q_{\varphi}(\spe{k},q)$.
\item  $Q_{\varphi}(\spe{h}, q^{-1}) = (-q)^{|I|}Q_{\varphi^{-1}}(\spe{h}, q)$.
\end{enumerate}
\end{theorem}
The last identity is a reciprocity result, and is due to \cite{stanley-poset} in the case of posets.

\section{Conclusion}

We conclude with questions:
\begin{enumerate}
\item Which properties of complexes are stable under the Hopf monoid operations in $\spe{\Phi}$? Do shellable complexes form a Hopf submonoid? What about Cohen-Macaulay complexes, or partitionable complexes?
\item What properties of a forbidden composition complex allow us to conclude that the triune quasisymmetric function is positive in the basis of fundamental quasisymmetric functions? This question is interesting: The complex $(\Gamma, \Delta)$ in Example \ref{ex:comp} part 2 has the feature that the Ehrhart quasisymmetric functions of $\Gamma$ and $\Delta$ are not $F$-positive, but the triune quasisymmetric function for $(\Gamma, \Delta)$ is $F$-positive.
\item If we linearize $\Phi$, what other natural bases does it possess? 
\end{enumerate}
Forbidden composition complexes, and triune quasisymmetric functions merit further study, as these geometric objects and their symmetric function invariants can be approached from three distinct perspectives.
\bibliographystyle{abbrvnat}
\bibliography{fpsac2016}

\end{document}